# Modelo matemático optimizado para la predicción y planificación de la asistencia sanitaria por la COVID-19

José Manuel Garrido[1,2,3,*], David Martínez-Rodríguez[4], Fernando Rodríguez-Serrano[1,2], José Miguel Pérez-Villares[5], Andrea Ferreiro-Marzal[3], María del Mar Jiménez-Quintana[5], Grupo de Estudio COVID-19_Granada [#], Rafael Jacinto Villanueva[4]

[1] Instituto de Investigación Biosanitaria ibs.GRANADA, Granada, España
[2] Instituto de Biopatología y Medicina Regenerativa (IBIMER), Universidad de Granada, Granada, España
[3] Servicio de Cirugía Cardiovascular, Hospital Virgen de las Nieves, Granada, España
[4] Instituto Universitario de Matemática Multidisciplinar, Universitat Politècnica de València, Valencia, España
[5] Servicio de Medicina Intensiva, Hospital Universitario Virgen de las Nieves, Granada, España

# Grupo de Estudio COVID-19_Granada: José Manuel Garrido Jiménez, Antonio Cansino Osuna, Juan Carlos Carrillo Santos, Cristina Carvajal Pedrosa, Mª Angeles García Rescalvo, Francisco José Guerrero García, Sebastián Manzanares Galán, Francisco Marti Jiménez, Manuel Enrique Reyes Nadal, Pedro Manuel Ruiz Lorenzo, José Luis Salcedo Lagullón, Indalecio Sánchez-Montesinos García, Fernando Rodríguez-Serrano

* Corresponding author. Departamento de Cirugía y sus Especialidades, Universidad de Granada, Granada, España, Avd. de la Investigación 11, 18016, Granada, España.
E-mail: josemgarrido@ugr.es

**Resumen**

**Objetivo**

La pandemia de la COVID-19 ha supuesto una amenaza de colapso de los servicios hospitalarios y de UCI, y una reducción de la dinámica asistencial de pacientes afectados por otras patologías. El objetivo fue desarrollar un modelo matemático diseñado para optimizar las predicciones relacionadas con las necesidades de hospitalización e ingresos en UCI por la COVID-19.

**Diseño**

Estudio prospectivo.

**Ámbito**

Provincia de Granada (España).

**Pacientes**

Pacientes consecutivos de COVID-19 hospitalizados, ingresados en UCI, recuperados y fallecidos desde el 15 de marzo hasta el 22 de septiembre de 2020.

**Intervenciones**

Desarrollo de un modelo matemático tipo SEIR capaz de predecir la evolución de la pandemia considerando las medidas de salud pública establecidas.

**Variables de interés**

Número de pacientes infectados por SARS-CoV-2, y hospitalizados e ingresados en UCI por la COVID-19.

**Resultados**

A partir de los datos registrados hemos podido desarrollar un modelo matemático que refleja el flujo de la población entre los diferentes grupos de interés en relación a la COVID-19. Esta herramienta nos ha permitido analizar diferentes escenarios basados en medidas de restricción socio-sanitarias, y pronosticar el número de infectados, hospitalizados e ingresados en UCI hasta el mes de mayo de 2021.

**Conclusiones**

El modelo matemático es capaz de proporcionar predicciones sobre la evolución de la COVID-19 con suficiente antelación como para poder conjugar los picos de prevalencia y de necesidades de asistencia hospitalaria y de UCI, con la aparición de ventanas temporales que posibiliten la atención de enfermos no-COVID.

# Mathematical model optimized for prediction and health care planning for COVID-19

## Abstract

### Objective

The COVID-19 pandemic has threatened to collapse hospital and ICU services, and it has affected the care programs for non-COVID patients. The objective was to develop a mathematical model designed to optimize predictions related to the need for hospitalization and ICU admission by COVID-19 patients.

### Design

Prospective study.

### Setting

Province of Granada (Spain)

### Population

Consecutive COVID-19 patients hospitalized, admitted to ICU, recovered and died from March 15 to September 22, 2020.

### Study variables

The number of patients infected with SARS-CoV-2 and hospitalized or admitted to ICU for COVID-19.

### Results

The data reported by hospitals was used to develop a mathematical model that reflects the flow of the population among the different interest groups in relation to COVID-19. This tool has allowed us to analyse different scenarios based on socio-health restriction measures, and to forecast the number of people infected, hospitalized and admitted to the ICU until May 2021.

### Conclusions

The mathematical model is capable of providing predictions on the evolution of the COVID-19 sufficiently in advance as to anticipate the peaks of prevalence and hospital and ICU care demands, and also the appearance of periods in which the care for non-COVID patients could be intensified.

**Introducción**

Los coronavirus causan enfermedades respiratorias e intestinales en numerosas especies animales. En humanos, cuatro coronavirus producen infecciones respiratorias de las vías altas (OC43, HKU1, 229E y NL63) y dos pueden causar síndromes respiratorios severos (SARSCoV-1 y MERS-CoV)[1]. Sin embargo, el pasado mes de diciembre las autoridades chinas informaron de diferentes casos de síndrome respiratorio en Wuhan que posteriormente fueron atribuidos a infecciones del nuevo coronarivus 2 del síndrome respiratorio agudo (SARS-CoV-2), el agente causante de la enfermedad por coronarivus 2019 (COVID-19)[2,3]. Desde la declaración de pandemia realizada por la Organización Mundial de la Salud, y hasta el 29 de noviembre de 2020, se han confirmado 61.869.330 casos acumulados y 1.448.896 fallecimientos en todo el mundo. En dicho período, en España se han registrado 1.628.208 casos y 44.668 fallecimientos[4], y todos los indicadores señalan a España como uno de los países más afectados por la enfermedad[5].

Hasta el momento, las principales líneas de tratamiento sugeridas contra el SARS-CoV-2 han incluido esteroides, plasma de convalecientes, anticuerpos neutralizantes, interferón, remdesivir, anticoagulantes, cloroquina/hidroxicloroquina, favipiravir y lopinavir/ritonavir, vitamina D, entre otros agentes. Sin embargo, existe una necesidad urgente de desarrollar nuevas terapias eficaces y especialmente vacunas que sean capaces de generar una suficiente inmunización comunitaria frente al SARSCoV-2[6-9]. Por ello, la aplicación de medidas no farmacéuticas, como distanciamiento social, el empleo de mascarillas faciales, la mejora de los hábitos de higiene, confinamientos perimetrados, confinamientos domiciliarios, cierre de servicios no esenciales, restricciones de movilidad, etc., han cobrado especial importancia por generar un impacto directo en la velocidad de propagación de la enfermedad[10-12]. De hecho, los indicadores de evolución de la COVID-19 mejoraron sensiblemente en España a las dos semanas de la declaración institucional de cuarentena del 14 de marzo de 2020, y mostraron además cómo en regiones que se encontraban en estadíos iniciales de pandemia en el momento del confinamiento, como Ceuta y Melilla, presentaron tasas de defunción muy inferiores respecto a otras que partían de una transmisión más acentuada, como Cataluña. Lo mencionado anteriormente refleja el importante efecto diferencial de las medidas cuando se adoptan de forma temprana[13]. Debido a la evolución de la pandemia a nivel de la provincia de Granada (Andalucía, España), se ha establecido un periodo de cierre de servicios no esenciales de dos semanas a partir del día 10 de noviembre de 2020[14].

Además de las consecuencias derivadas de la afectación poblacional por la COVID-19, la pandemia ha reducido de forma muy importante la dinámica asistencial de pacientes aquejados de otras patologías. Las consecuencias en términos de morbi-mortalidad se podrán ir cuantificando a lo largo del tiempo comparando con series retrospectivas pre-pandémicas. En este sentido, el Grupo de Trabajo de Planificación, Organización y Gestión de la SEMICYUC ha recomendado la implantación de un plan de contingencia en el servicio de medicina intensiva (SMI) integrado en el plan de contingencia local. Los profesionales del SMI deben ser reconocidos como miembros activos en los comités locales encargados de adaptar la planificación hospitalaria a la situación epidemiológica, garantizando una adecuada gestión de los recursos disponibles y de los necesarios[15]. Los planes deben permitir preservar programas con relevancia clínica y estratégica, como son los programas de trasplantes o los cuidados intensivos postoperatorio para Cirugía Cardiovascular, de manera adaptada al escenario epidemiológico actual y teniendo como objetivo final la recuperación de niveles de actividad similares a los alcanzados en nuestro país en situación pre-pandémica.

Para conseguir una adecuada planificación es necesario emplear herramientas que permitan predecir la evolución de la COVID-19 en función de la situación de partida, y de las medidas no farmacéuticas y de salud pública instauradas en cada momento, y que sean capaces de alertar sobre los posibles escenarios de transmisión. Las herramientas deben permitir el diseño de la aplicación y temporalización de medidas con suficiente antelación como para poder conjugar los picos de prevalencia y de necesidades de asistencia hospitalaria y de UCI por la COVID-19, con la aparición de ventanas temporales que posibiliten la atención de enfermos no-COVID. La SEMICYUC ha realizado una planificación de posibles escenarios utilizando el software FluSurge 2.0[16], que ha sido desarrollado por "Centers for Disease Control and Prevention (CDC)". La aplicación permite realizar cálculos aproximados de la demanda de servicios en una situación de pandemia moderada y grave, considerando la población en riesgo, los recursos hospitalarios disponibles y los supuestos de curso epidemiológico del proceso pandémico. Sin embargo, FluSurge 2.0 se ha diseñado específicamente para valorar el efecto de una pandemia por el virus de la influenza y solo se ha validado para dicho ámbito. Además, no permite rediseñar el flujo de pacientes a nivel hospitalario ni incorporar nuevos factores que afecten a la dinámica de transmisión, como es una potencial campaña de vacunación.

En el presente artículo presentamos un modelo matemático capaz de predecir la dinámica de transmisión de la COVID-19, que ha sido calibrado y validado empleando datos proporcionados por los hospitales de la provincia de Granada, los cuales dan conjuntamente una cobertura asistencial a 914.678 habitantes[17]. Además, presentamos estimaciones de diferentes escenarios basados en las medidas de contención que se están activando en Andalucía, y en otras localizaciones españolas, lo que nos ofrece la posibilidad de diseñar estrategias que permitan prever la utilidad de las medidas, y los picos de infectados, hospitalizados e ingresos en UCI. Todo ello, nos permitiría diluir a lo largo del tiempo la prevalencia de la enfermedad evitando la saturación de los servicios hospitalarios, especialmente de UCIs, y prever periodos de valle, o menor carga COVID, para programar la atención de enfermos no-COVID reduciendo así el coste de oportunidad. De igual modo, el diseño del calendario de aplicación de las medidas de contención socio-sanitaria que pueden realizarse debe contemplar también su impacto durante los periodos de mayor actividad económica, como Navidad o Semana Santa.

**Pacientes y métodos**

*Pacientes*

Para la calibración del modelo, hemos empleado los datos registrados por los hospitales de Granada H.U. Virgen de las Nieves, H.U. Clínico San Cecilio, H.U. Santa Ana en Motril, H.U. de Baza, H. de San Rafael y H.L.A. Inmaculada. Para ello, contamos con la aprobación por parte del Comité de Ética de la Investigación Biomédica de la provincia de Granada dependiente de la Consería de Salud y Familias de la Junta de Andalucía. Recopilamos el número de hospitalizados, ingresados en UCI, recuperados y fallecidos desde el 15 de marzo hasta el 22 de septiembre de 2020.

*Modelo matemático, calibración y validación*

Hemos implementado un modelo SEIR (susceptible, expuesto, infectado y recuperado) diseñado específicamente para describir la dinámica de la epidemia a nivel poblacional y a nivel de circuito hospitalario en relación a los pacientes de COVID-19 (hospitalizaciones en planta y UCI), ya que es el aspecto más limitante a la hora de hacer frente a la pandemia, debido a los recursos materiales y personales que requiere (**Figura 1**). La tabla 1 recoge los diferentes grupos en los que se puede segregar a la población respecto a la infección y el circuito hospitalario, junto con las ecuaciones en diferencias que describen la dinámica de cada grupo a lo largo del tiempo. Cada individuo puede pertenecer al grupo susceptible (S), cuarentena (Q), latente (L), infeccioso (I), recuperado (R), hospitalizado (H), ingresado en UCI (U), fallecido (F), hospitalizado tras pasar por la UCI (HU) o alta (A). El paso entre grupos viene determinado por las tasas de transición $q_s$, $s_q$, $l_i$, $i_r$, $i_h$, $i_u$, $h_u$, $h_f$, $h_a$, $u_f$, $u_{hu}$, $h_f$ y $hu_a$. $\beta$ es la tasa de transmisión entre S e I, y su valor es proporcional a la magnitud del número de reproducción básico $R_0$, según la expresión: $R_0 = \beta / (i_r + i_h + i_u)$. Durante el proceso de calibración del modelo, a partir de los registros hospitalarios pudimos determinar el valor de las diferentes tasas de transición y de la tasa de transmisión ($\beta$) necesarias para que el modelo sea capaz de describir la situación específica de la provincia de Granada, empleando para ello el algoritmo de optimización Novelty Swarm implementado en Pyhton3[18]. El valor de $s_q$ y $q_s$ se estableció atendiendo al cambio en los flujos de desplazamientos recogidos en el "Informe de Movilidad Local sobre la COVID-19" que proporciona la compañía Google[19], en el periodo de confinamiento domiciliario decretado por el gobierno español el 14 de marzo de 2020[20]. La calibración se realizó empleando el número real de registros de hospitalizados, ingresados en UCI, recuperados y fallecidos en los hospitales de la provincia de Granada, entre el 15 de marzo y el 22 de septiembre de 2020. Una vez finalizada la calibración, validamos el modelo comparando los datos predichos por el modelo y los datos reales registrados entre el 23 de septiembre y 7 de noviembre de 2020.

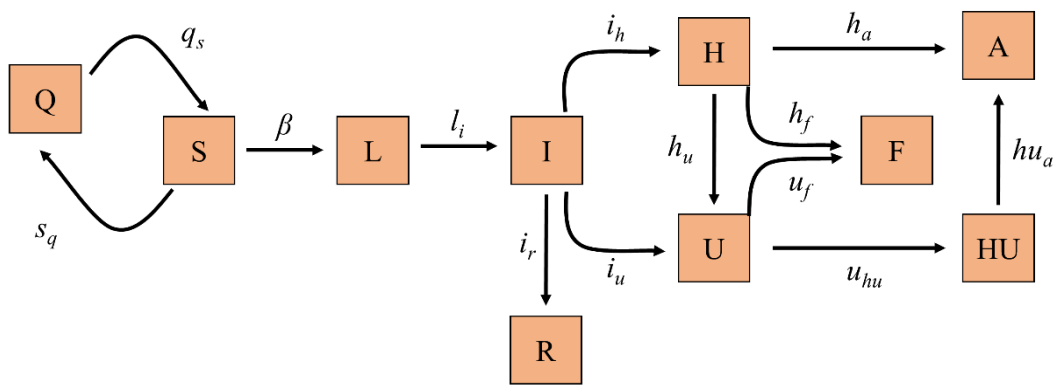

**Figura 1**. Diagrama del modelo SEIR implementado para la predicción de la transmisión de la COVID-19. Cada caja representa uno de los diferentes grupos en los que se puede segregar a la población respecto a la infección y la enfermedad. Los parámetros de las flechas representan las tasas de transmisión entre los grupos anexos. Susceptible (S); Cuarentena (Q); Latente (L); Infeccioso (I); Recuperado (R); Hospitalizado (H); UCI (U); Fallecido (F); Hospitalizado tras UCI (HU); Alta (A).

**Tabla 1.** Grupos de población en relación a la infección por el virus SARS-CoV-2 y a la evolución de la COVID-19, y ecuaciones que predicen en cada momento la cuantificación de cada grupo.

| Grupo | Descripción | Ecuación* |
|---|---|---|
| Susceptible (S) | Sin exposición al SARS-CoV-2 | $S(t+1) = S(t) + q_s(t) - s_q(t) - \beta(t)S(t)(I(t)/P_T)$ |
| Cuarentena (Q) | En confinamiento domiciliario | $Q(t+1) = Q(t) + s_q(t) - q_s(t)$ |
| Latente (L) | Infectado de SARS-CoV-2 pero no infeccioso | $L(t+1) = L(t) + \beta(t)S(t)(I(t)/P_T) - l_i L(t)$ |
| Infeccioso (I) | Infectado de SARS-CoV-2 e infeccioso | $I(t+1) = I(t) + l_i L(t) - (i_r(t) + i_h(t) + i_u(t))I(t)$ |
| Recuperado (R) | Paciente no hospitalizado que supera la COVID-19 | $R(t+1) = R(t) + i_r(t)I(t)$ |
| Hospitalizado (H) | Paciente hospitalizado en planta | $H(t+1) = H(t) + i_h(t)I(t) - (h_u(t) + h_f(t) + h_a(t))H(t)$ |
| UCI (U) | Paciente en UCI | $U(t+1) = U(t) + i_u(t)I(t) + h_u(t)H(t) - (u_f(t) + u_{hu}(t))U(t)$ |
| Fallecido (F) | Paciente fallecido a causa de la COVID-19 | $F(t+1) = F(t) + h_f(t)H(t) + u_f(t)U(t)$ |
| Hospitalizado tras UCI (HU) | Paciente que pasa de UCI a otras dependecias hospitalarias por mejoría | $HU(t+1) = HU(t) + u_{hu}(t)U(t) - hu_a(t)HU(t)$ |
| Alta (A) | Paciente dado de alta del hospital | $A(t+1) = A(t) + h_a(t)H(t) + hu_a(t)HU(t)$ |

* $q_s$, $s_q$, $l_i$, $i_r$, $i_h$, $i_u$, $h_u$, $h_f$, $h_a$, $u_f$, $u_{hu}$, $h_f$ y $hu_a$ se corresponden con las tasas de transición entre grupos de sujetos/pacientes. $\beta$ es la tasa de transmisión entre S e I, y su valor es proporcional a la magnitud del número básico de reproducción $R_0$, según la expresión: $R_0 = \beta / (i_r + i_h + i_u)$. $P_T$ es la población de la provincia de Granada.

*Predicciones proporcionadas por el modelo*

Hemos generado 3 escenarios atendiendo a medidas no farmacéuticas para analizar la evolución más probable de la pandemia, y para determinar las condiciones más favorables para conjugar los picos y valles de prevalencia y de necesidades de hospitalización y de ingresos en UCIs. Esto serviría tanto para evitar la saturación de los recursos sanitarios, como para favorecer la asistencia de los enfermos no-COVID, posibilitando además periodos de relajación de las medidas para promover la actividad económica.

El escenario inicial representa la evolución predicha atendiendo a las restricciones establecidas para la provincia de Granada, que se encuentra en nivel 4 fase 2, lo que implica el cese temporal de servicios no esenciales, restricción del horario comercial y de movilidad, y confinamiento perimetral[14]. Por el momento, estas medidas han sido decretadas para dos semanas a partir del día 10 de noviembre de 2020. El efecto de dichas medidas se podrá cuantificar a corto plazo, pero para nuestras simulaciones hemos considerado que tendrán un impacto que, en el mejor de los casos, podría igualar al valor de $R_t$ que aparece actualmente en Cataluña, tras la aplicación de medidas similares desde el 30 de octubre de 2020, y cuyo valor se situó en torno a 0,8[21]. Los otros dos escenarios adicionales han sido escogidos entre muchas simulaciones ya que permiten prever los efectos de la dilatación temporal de dichas medidas y el establecimiento de diferentes periodos de restricciones. El modelo cuantifica para cada escenario el valor esperado diario en cada grupo (S, Q, L, I, R, H, U, F, HU y A), y el intervalo de confianza del 95% (percentiles 2.5 – 97.5).

**Resultados**

El modelo desarrollado nos ha permitido establecer diferentes escenarios de aplicación de medidas de restricción y prever la evolución del número de infectados, ingresos hospitalarios y en UCI (**Tabla 2**), considerando una $R_0$ de 0,8. La fase de calibración y validación del modelo ha demostrado la validez de las predicciones proporcionadas por el modelo tras comparar los casos esperados y registrados de hospitalizaciones e ingresos en UCI en la provincia de Granada durante el periodo 23 de septiembre y 7 de noviembre de 2020 (**Figura 2**). En la Figura 2 podemos observar que, a pesar de que no todos los puntos se encuentran dentro del intervalo de confianza, el crecimiento de ambas curvas tiene una forma muy similar. Es por ello por lo que el modelo puede ser capaz de obtener una tendencia clara de la evolución del circuito hospitalario que permita valorar de forma cualitativa la evolución de la pandemia. A continuación, se presentará un resumen de la evolución prevista en los tres escenarios modelizados, y todas las referencias numéricas representan el resultado de la media ofrecida por el modelo.

El escenario 1 se corresponde con la situación de restricciones de restricción socio-sanitaria establecidas durante dos semanas a partir del día 10-11-2020. En este contexto, el número de infectados alcanzaría un pico de 58.379 afectados el 06-02-21, un dato de prevalencia que es más de 3 veces superior respecto a lo registrado durante el pico que apareció en marzo (18.448 el 20-03-20). Sólo encontraríamos una breve reducción de la tendencia de aumento de infectados a finales de octubre (26-11-20), antes de volver a experimentar un incremento exponencial de casos. En el caso del número de hospitalizaciones e ingresos en UCI, la medida de restricción de 2 semanas no conseguiría generar un valle que reduzca la presión hospitalaria, y encontraríamos un pico de 1.946 (fecha 20-02-21) y 310 (fecha 25-02-21) casos, respectivamente. Dicho volumen de casos superaría de forma importante el número de pacientes hospitalizados y en UCI atendidos en el pasado periodo marzo-abril de 2020 (**Figura 3**). Cabe matizar que la presión hospitalaria durante la primea ola, aunque elevada, no llegó a saturar los recursos hospitalarios de la provincia, ya que se decretó confinamiento domiciliario cuando la prevalencia en la provincia se encontraba en sus estadíos iniciales (22).

Si consideramos el escenario 2, con una dilatación de las restricciones hasta un total de 4 semanas desde el 10-11-2020, una situación que sería similar a lo acontecido en Cataluña hasta el momento, encontraríamos que el número de afectados se reduciría durante el periodo 13-11-20 al 10-12-20, pasando de 25.233 a 16.359 casos. No obstante, encontraríamos un pico elevado de prevalencia de 50.485 casos más tarde respecto al pico estimado para el escenario anterior (02-03-21). En el caso de hospitalizados e ingresos en UCI encontraríamos un comportamiento parecido, con una reducción acentuada en el número de casos en los períodos 26-11-20 / 17-12-20, pasando de 733 (IC95% 509-1.048) a 641 (IC95% 449-900) casos, y 01-12-20 / 18-12-20, pasando de 112 (IC95% 75-158) a 106 (IC95% 71-148) casos, respectivamente. En ambos grupos encontraríamos un pico en marzo de 2021, que superaría más de 4 veces el número de pacientes hospitalizados y en UCI atendidos en el pasado periodo marzo-abril de 2020 (**Figura 4**).

El escenario 3 presenta la evolución de casos antendiendo a un calendario de aplicación de medidas de restricción socio-sanitaria por fases, después de analizar numerosos escenarios posibles. Este contexto contemplaría 3 periodos: 4 semanas a partir del día 10/11/20 de restricciones en servicios no esenciales; 4 semanas a partir del día 11/01/21 con un confinamiento poblacional del 70%; y 2 semanas a partir del día 15/03/21 de restricciones en servicios no esenciales. Este escenario predice una evolución creciente y más suave que los escenarios anteriores hasta enero, momento en el que encontraríamos 30.429 infectados (13-01-21), 889 hospitalizados (19-01-21) y 137 ingresos en UCI (21-01-21). Además, en diciembre encontraríamos una reducción relevante de casos en los tres grupos, que estaría más acentuada en el grupo de infectados, que pasaría de 25.233 casos el 13-11-20 a 15.753 el 12-12-20. A partir de enero, el número de casos se iría reduciendo con puntos alternantes de picos y valles hasta llegar a mayo, cuándo encontraríamos 10.922 infectados (07-05-21), 324 hospitalizados (01-05-21) y 52 ingresos en UCI (01-05-21) (**Figura 5**).

**Tabla 2**. Número de infectados de SARS-CoV-2, hospitalizados e ingresados en UVI previstos por el modelo matemático para seis escenarios diferenciados en el calendario y en la duración de la aplicación de medidas de restricción para la provincia de Granada.

|  | (Escenario 1) 10/11/20 (2sem.) | | (Escenario 2) 10/11/20 (4 sem.) | | (Escenario 3) 10/11/20 (4 sem.) 11/01/21 (4 sem.) 15/03/21 (2 sem.) | |
|---|---|---|---|---|---|---|
| **Infectados** | 20-03-20 | 18.448 (17.620-19.331) | 20-03-20 | 18.449 (17.620-19.331) | 20-03-20 | 18.449 (17.620-19.331) |
| | 06-06-20 | 625 (590-664) | 06-06-20 | 625 (590-664) | 06-06-20 | 625 (590-664) |
| | 13-11-20 | 25.206 (16.907-36.114) | 13-11-20 | 25.233 (16.908-36.109) | 13-11-20 | 25.233 (16.908-36.109) |
| | 26-11-20 | 21.191 (14.390-29.963) | 10-12-20 | 16.359 (11.379-22.504) | 12-12-20 | 15.753 (10.999-21.577) |
| | 06-02-21 | 58.379 (45.942-72.245) | 02-03-21 | 50.485 (40.249-61.038) | 13-01-21 | 30.429 (20.582-42.286) |
| | | | | | 12-02-21 | 6.280 (4.528-8.079) |
| | | | | | 18-03-21 | 10.669 (7.871-13.230) |
| | | | | | 04-04-21 | 7.640 (5.834-9.243) |
| | | | | | 07-05-21 | 10.922 (8.896-12.923) |
| **Hospitalizados** | 31-03-20 | 399 (387-413) | 31-03-20 | 399 (388-413) | 31-03-20 | 399 (388-413) |
| | 30-06-20 | 13 (10-16) | 30-06-20 | 13 (10-16) | 30-06-20 | 13 (10-16) |
| | 20-02-21 | 1.946 (509-2.418) | 26-11-20 | 733 (509-1.048) | 26-11-20 | 733 (509-1.048) |
| | | | 17-12-20 | 641 (449-900) | 20-12-20 | 623 (437-872) |
| | | | 16-03-21 | 1.691 (1.366-2.052) | 19-01-21 | 889 (592-1.226) |
| | | | | | 07-03-21 | 314 (224-406) |
| | | | | | 24-03-21 | 340 (248-434) |
| | | | | | 15-04-21 | 296 (224-367) |
| | | | | | 01-05-21 | 324 (247-394) |
| **UCI** | 05-04-20 | 56 (51-62) | 05-04-20 | 56 (51-62) | 05-04-20 | 56 (51-62) |
| | 25-07-20 | 2 (2-3) | 25-07-20 | 2 (2-3) | 25-07-20 | 2 (2-3) |
| | 25-02-21 | 310 (248-387) | 01-12-20 | 112 (75-158) | 01-12-20 | 112 (75-158) |
| | | | 18-12-20 | 106 (71-148) | 21-12-20 | 103 (69-144) |
| | | | 21-03-21 | 270 (217-332) | 21-01-21 | 137 (89-192) |
| | | | | | 15-03-21 | 56 (39-72) |
| | | | | | 23-03-21 | 56 (41-72) |
| | | | | | 18-04-21 | 50 (37-61) |
| | | | | | 01-05-21 | 52 (39-63) |

* Para cada uno de los 3 escenarios se presenta la fecha de implantación y la duración de las medidas de restricción, así como la fecha y el número de pacientes y el intervalo de confianza de infectados, hospitalizados o en UCI correspondientes a los puntos de inflexión, pico (máximos) o valle (mínimos), que aparecerían a lo largo de la predicción. Sem, semanas.

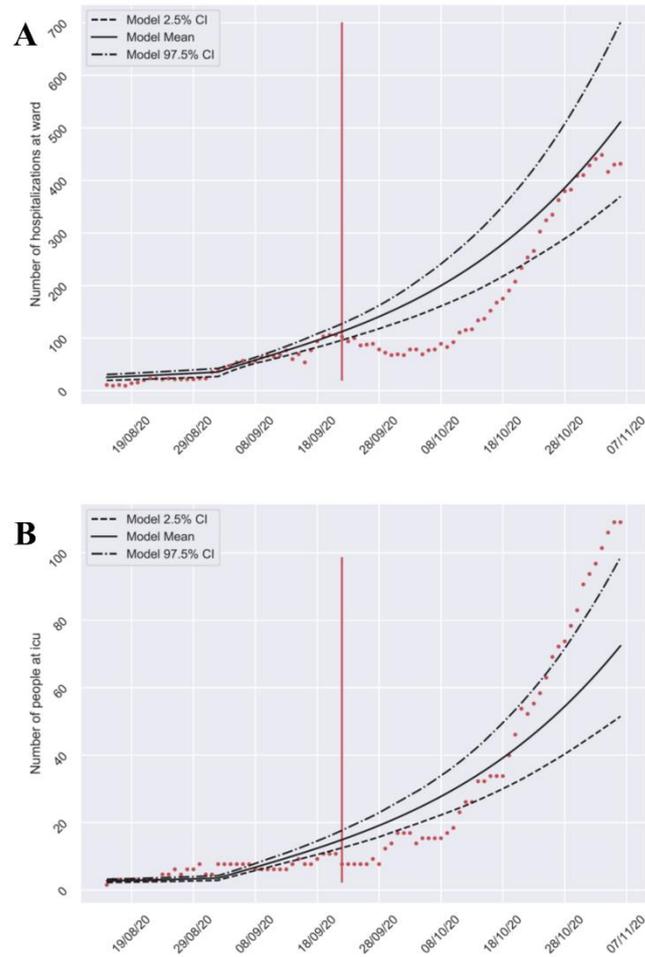

**Figura 2**. Validación del modelo matemático comparando el número de hospitalizados (A) e ingresados en UCI (B) en los hospitales de la provincia de Granada durante el periodo 23 de septiembre y 7 de noviembre de 2020 (puntos), respecto a la media e intervalo de confianza predicho por el modelo para dicho periodo.

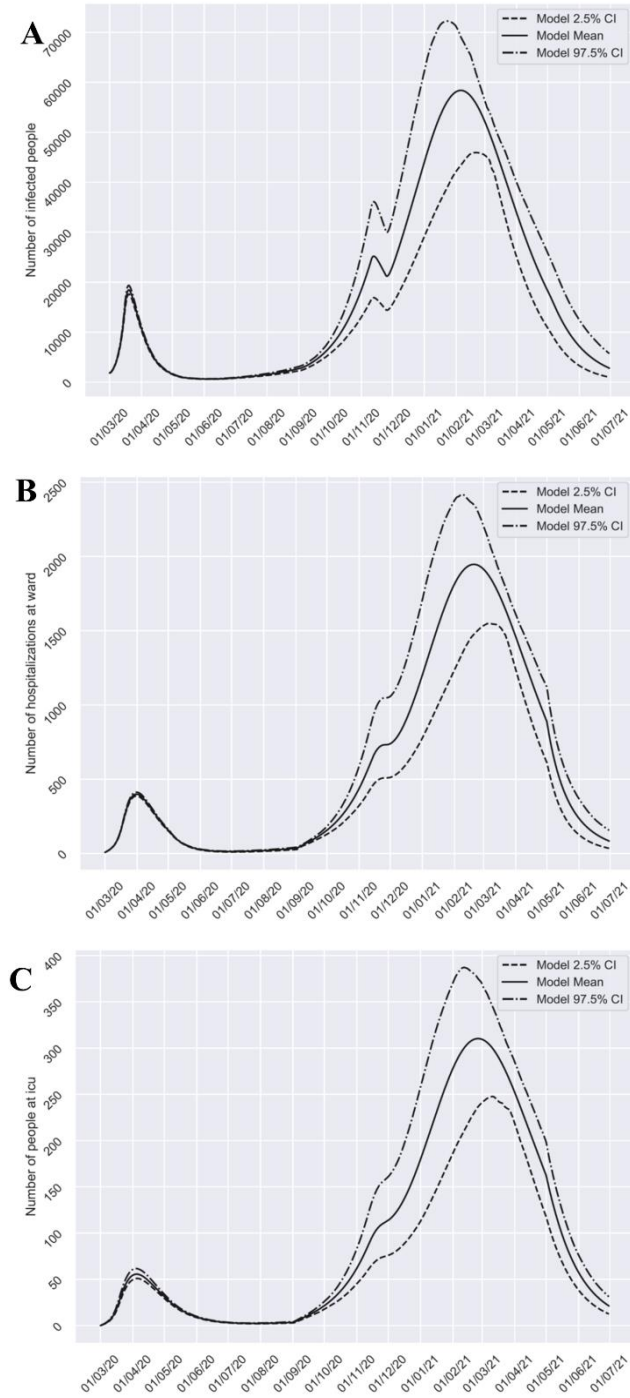

**Figura 3**. Predicciones de prevalencia de infectados por SARS-CoV-2, hospitalizados e ingresos en UCI para la provincia de Granada considerando 2 semanas de restricciones en servicios no esenciales desde el 10-11-2020 (Escenario 1). Se representa el número de infectados (A), hospitalizados (B) e ingresos en UCI (C) a lo largo del tiempo.

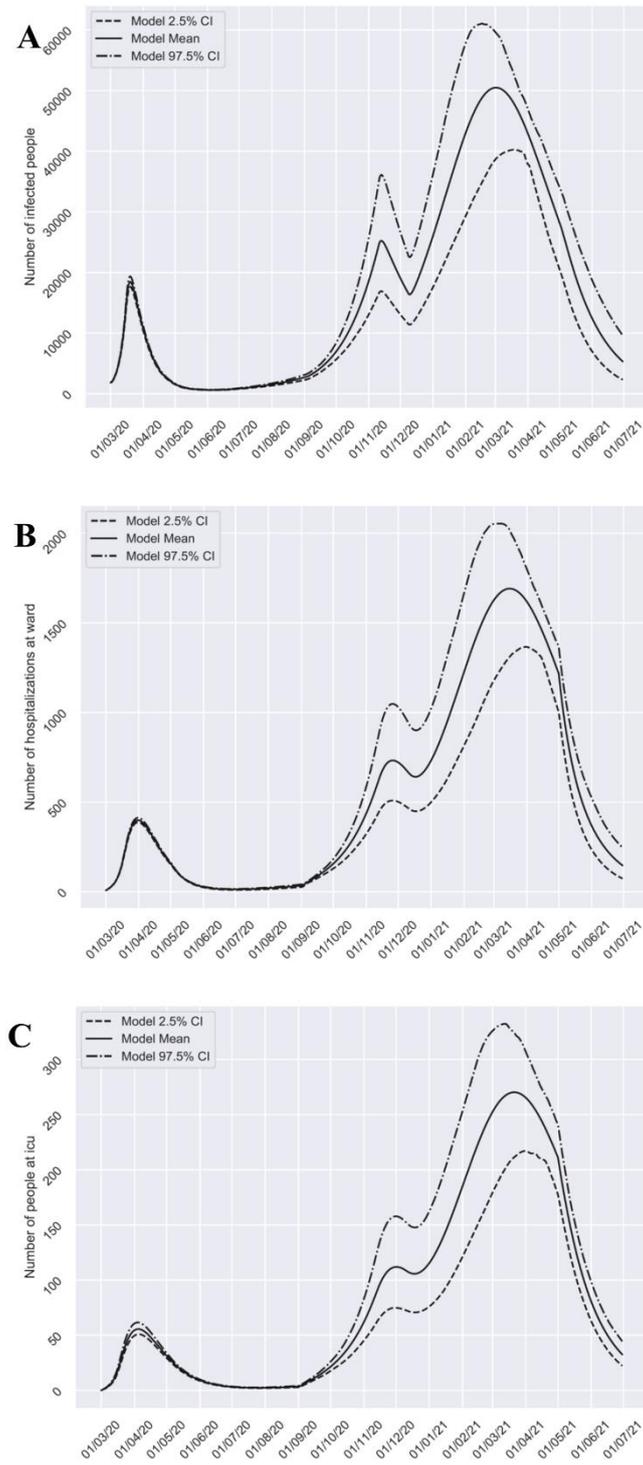

**Figura 4**. Predicciones de prevalencia de infectados por SARS-CoV-2, hospitalizados e ingresos en UCI para la provincia de Granada considerando 4 semanas de restricciones en servicios no esenciales desde el 10-11-2020 (Escenario 2). Se representa el número de infectados (A), hospitalizados (B) e ingresos en UCI (C) a lo largo del tiempo.

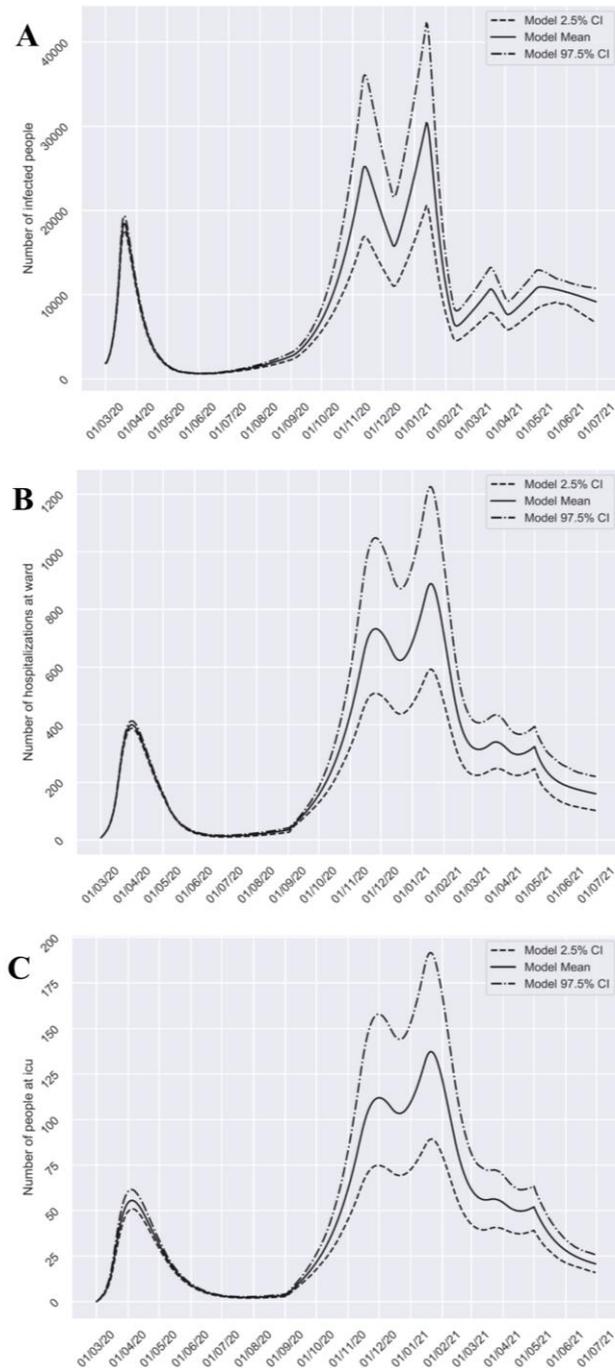

**Figura 5**. Predicciones de prevalencia de infectados por SARS-CoV-2, hospitalizados e ingresos en UCI para la provincia de Granada considerando: 4 semanas a partir del día 10/11/20 de restricciones en servicios no esenciales; 4 semanas a partir del día 11/01/21 con un confinamiento poblacional del 70%; y 2 semanas a partir del día 15/03/21 de restricciones en servicios no esenciales (Escenario 3). Se representa el número de infectados (A), hospitalizados (B) e ingresos en UCI (C) a lo largo del tiempo.

**Discusión**

La compleja situación asistencial debido a la actual pandemia COVID-19 se potencia por la importante afectación del normal funcionamiento de la atención sanitaria y hospitalaria[23-25], y por el elevado coste de oportunidad que se ha producido en relación a muchas patologías graves que han dejado de atenderse de forma óptima. El adecuado tratamiento y seguimiento de ciertos grupos de pacientes cobra especial importancia debido a que parecen representar comorbilidades correlacionadas con las necesidades de ingresos hospitalarios y en UCI de infectados por SARS-CoV2, tales como la hipertensión, enfermedades crónicas cardíacas, diabetes, enfermedades crónicas pulmonares y obesidad[26].

El control de la primera ola de la COVID-19 se realizó mediante la aplicación de medidas extremas de restricción poblacional al amparo del decreto de estado de alarma del 14 de marzo de 2020. Las medidas establecidas incluyeron un confinamiento domiciliario de 7 semanas, seguido de cuatro fases de desconfinamiento hasta alcanzar la llamada "nueva normalidad" el 21 de junio de 2020[27]. Durante este periodo se produjo una reorganización de los Sistemas Sanitarios para encajar la atención prioritaria de pacientes COVID-19, sobrepasando ampliamente los límites estratégicos y operacionales de dichos sistemas, y relegando a un segundo plano la atención del resto de patologías graves y no graves que representan la habitual cartera de servicios ofertada por nuestros centros sanitarios. La desescalada progresiva de restricciones no ha llegado en ningún momento a la recuperación de un contexto social similar al previamente existente. Por el contrario, la obligatoriedad de medidas higiénico-sanitarias ha ido oscilando entre niveles variables de intensidad.

La presencia del virus SARS-CoV-2 en nuestro medio se ha mantenido constante tras la finalización de la primera ola epidémica en España. Como consecuencia de ello, la endemización de la enfermedad en Europa es un hecho no discutible, alternándose fases temporales de mínima incidencia de la infección, con otras fases de estallido epidémico exponencial[28]. Este patrón epidemiológico que describe la dinámica temporal de la pandemia COVID-19 se mantendrá, con elevada probabilidad, hasta la aparición de una vacuna efectiva y/o el desarrollo de una apropiada inmunidad poblacional[29].

Asimismo, y debido a la heterogénea incidencia de la infección y a los diferentes niveles de colapso/puesta al límite de los sistemas sanitarios en las distintas regiones de España, la elaboración de una estrategia única para afrontar el impacto socio-sanitario del virus durante los próximos años resulta especialmente compleja. Sin embargo, y a pesar de la enorme dificultad que representa el disponer de una estrategia conjunta, es necesario resaltar dos aspectos fundamentales que deben tenerse en cuenta en la elaboración de cualquier planificación concreta. Por un lado, hay que desarrollar una vía asistencial COVID-19 que asegure una adecuada atención de estos pacientes, tanto en régimen ambulatorio como de ingreso hospitalario o en UCI, y por otro lado, el desarrollo de una vía asistencial no COVID-19, que nos permita reducir al máximo el coste de oportunidad relativo al resto de patologías graves y que consiga amortiguar el incremento sustancial de la morbi-mortalidad asociada.

Sobre esta base, y tras la finalización de la fase de confinamiento poblacional en el pasado mes de junio, los gestores sanitarios han intentado fortalecer y reorganizar los recursos sanitarios existentes, con objeto de satisfacer los dos aspectos fundamentales antes señalados. A pesar de ello, durante la llamada segunda oleada de la pandemia COVID-19, el incremento exponencial de pacientes con necesidad de ingreso hospitalario y en UCI ha vuelto a condicionar que la mayor parte de los recursos sanitarios que se dirigen hacia el tratamiento de los pacientes COVID-19. El coste de oportunidad acumulado del resto de patologías graves durante el tiempo que duró la primera oleada COVID-19 y la fase de recuperación/normalización de las instituciones sanitarias puso de manifiesto la imposibilidad de mantener una adecuada atención poblacional fuera del desbordante COVID-19. Esta situación de incapacidad de mantener ambos circuitos (COVID-19/no-COVID), incluso con el refuerzo de las estructuras sanitarias es un hecho evidenciable por el porcentaje de recursos dirigidos a pacientes COVID-19[30]. Además, el problema ético, derivado de la polarización del sistema hacia el tratamiento de la pandemia, no ha sido resuelto.

La identificación del comportamiento poblacional del virus en picos y valles resulta de especial trascendencia para una correcta planificación. De forma lógica, y preservando una atención de calidad para los pacientes urgentes/emergentes de las distintas patologías graves, puede organizarse la asistencia como periodos alternantes de focalización asistencial. Así, durante los picos se intensificaría la atención COVID-19, reduciendo, proporcionalmente a la magnitud de los mismos, la asistencia sanitaria de los pacientes estables de otras patologías. Por el contrario, durante las fases de valle, intensificaríamos, por encima del estándar habitual, la actividad clínica dirigida a patologías no-COVID, aprovechando plenamente la ventana de oportunidad generada.

La problemática de este comportamiento se centra en conocer qué medidas de restricción poblacional son necesarias para alcanzar los valles, cuál es la duración efectiva de esos periodos de menor presión asistencial COVID-19, cuál es el punto umbral a partir del cual se progresa exponencialmente hacia un nuevo pico o cuál será la duración e intensidad del mismo. Además, hay que tener en cuenta que las distintas estrategias poblacionales y sanitarias que pueden adoptarse, la incertidumbre, junto al perjudicial impacto económico, aumentan la dificultad para el diseño y la aplicación de medidas efectivas.

En este sentido, el modelo matemático que presentamos en este artículo puede contribuir de manera importante en la toma de decisiones relacionadas con la aplicación de medidas y su calendario. Posibilita la planificación analizando los escenarios que las distintas estrategias de control puedan generar. El modelo ha sido diseñado para optimizar las predicciones relacionadas con las necesidades de hospitalización y de ingresos en UCI, aspectos más limitantes para mantener una adecuada atención de los circuitos COVID-19 y no-COVID. Además, el modelo tiene una construcción modular, de manera que es posible incorporar otros grupos en función de la aparición de nuevos factores que influyan de forma importante en la dinámica de transmisión del virus, como podría ser la aparición de una vacuna frente al SARS-CoV-2, y su porcentaje de eficacia asociado.

El escenario 1 predice la evolución de la pandemia considerando dos semanas de restricciones socio-sanitarias en servicios no esenciales. Podemos observar cómo resultaría insuficiente para cambiar la proyección de ingresos hospitalarios e ingresos en UCI de forma significativa y relevante a corto y medio plazo. Si bien es cierto que enlentece transitoriamente la curva, no genera impacto poblacional suficiente para cambiar la expectativa de saturación del sistema.

El escenario 2, que predice la evolución a partir de 4 semanas de restricciones, genera una respuesta positiva. Podemos observar que los efectos comienzan a tener un impacto significativo en diciembre. Así, se produciría un importante freno inicial del incremento exponencial de enfermos hospitalizados que se venía registrando durante los meses de octubre e inicio de noviembre. Sin embargo, encontraríamos un número de infectados, hospitalizados e ingresos en UCI exponencialmente ascendente hasta alcanzar un pico no asumible a lo largo de marzo de 2021. No obstante, este escenario genera un primer valle operativo que se extiende desde la última semana de noviembre a la última semana de diciembre, de 1 mes de duración, cuándo se produciría una reducción de presión asistencial, que permitiría redirigir parcialmente los recursos sanitarios hacia enfermos no-COVID.

El escenario 3 modeliza el efecto de 4 semanas a partir del día 10/11/20 de restricciones en servicios no esenciales, 4 semanas a partir del día 11/01/21 con un confinamiento poblacional del 70%, y 2 semanas a partir del día 15/03/21 de restricciones en servicios no esenciales. Lo analizado en el escenario 2, que representa la primera fase de este escenario, es de total aplicación para la evaluación del primer valle y de su potencial efecto beneficioso, concretado principalmente en el freno del crecimiento exponencial de la pandemia. Posteriormente, y después de alcanzar unas cifras de contagiados, pacientes hospitalizados y pacientes ingresados en UCI similares al periodo previo al del establecimiento de las medidas, se aplicaría una nueva fase de restricciones de 4 semanas (confinamiento poblacional del 70%) que generaría efectos claramente beneficiosos después de alcanzar un segundo pico de presión asistencial COVID-19, que se registraría en enero de 2021. A partir de ese momento, el impacto positivo de las restricciones comenzaría a ser pleno y se produciría una reducción progresiva del contagio poblacional con una disminución mucho muy acentuada y rápida del número de pacientes hospitalizados y de pacientes ingresados en UCI. Tal y como se puede observar en los datos proporcionados, este escenario generaría un muy importante periodo de presión asistencial decreciente COVID-19 desde enero, que llegaría a estabilizarse, en términos generales, a partir de marzo. Posteriormente, desde el mes de mayo entraríamos presumiblemente en la fase estimada de cambio estival de patrón de contagio y virulencia. Por ello, el escenario 3 combina distintas fases de restricciones poblacionales y de retirada de las mimas, consiguiendo un evidente efecto positivo y sostenible para el sistema sanitario. Se reduciría ostensiblemente el número de contagios poblacionales, de pacientes hospitalizados y de pacientes ingresados en UCI. Además, con la aplicación de dichas medidas, la curva predicha exponencial se desplaza hacia una modelización en picos y valles de presión asistencial, con una progresiva reducción de la carga de enfermedad COVID-19.

El modelo de previsión epidemiológica que planteamos nos permite evaluar el impacto de las distintas estrategias de restricción poblacional frente a la COVID-19, considerando su duración, intensidad y el contexto basal de incidencia y prevalencia de la COVID-19, así como prever el nivel de presión asistencial relativo al número de pacientes hospitalizados e ingresados en UCI. Todo ello, convierte nuestro modelo en una herramienta muy adecuada para diseñar planes de actuación y planes asistenciales a medio plazo, pudiendo incorporar cronológicamente los eventos de contención para determinar la duración e importancia de los valles de menor presión asistencial COVID-19, que podrían ser aprovechados para la asistencia de pacientes no-COVID. Además, nuestro modelo puede ser adaptado a otros núcleos poblaciones realizando una nueva calibración del modelo a partir de los correspondientes datos demográficos y de evolución local de la pandemia.

**Conflictos de intereses**

Ninguno

**Referencias**